\documentclass[oneside,10pt,]{article}
\usepackage{geometry}
\usepackage{ifthen}
\usepackage{etoolbox}
\usepackage{ifxetex,ifluatex}
\usepackage{graphicx}
\PassOptionsToPackage{usenames,dvipsnames,svgnames,table}{xcolor}
\usepackage{xcolor}
\usepackage{tcolorbox}
\tcbuselibrary{skins}
\tcbuselibrary{breakable}
\tcbuselibrary{raster}
\tcbset{ bwminimalstyle/.style={size=minimal, boxrule=-0.3pt, frame empty,
colback=white, colbacktitle=white, coltitle=black, opacityfill=0.0} }
\tcbset{ runintitlestyle/.style={fonttitle=\normalfont\bfseries, attach title to upper} }
\tcbset{ exercisespacingstyle/.style={before skip={1.5ex plus 0.5ex}} }
\tcbset{ blockspacingstyle/.style={before skip={2.0ex plus 0.5ex}} }
\tcbuselibrary{xparse}
\IfFileExists{latexrelease.sty}{}{\usepackage{fixltx2e}}
\ifthenelse{\boolean{xetex} \or \boolean{luatex}}{%
\ifxetex\usepackage{xltxtra}\fi
\ifluatex\usepackage{realscripts}\fi
\usepackage{fontspec}
\usepackage{polyglossia}
\setmainlanguage[variant=usmax]{english}
}{%
\usepackage[utf8]{inputenc}
\usepackage{lmodern}
\usepackage[T1]{fontenc}
}
\usepackage{amsmath}
\usepackage{amssymb}
\allowdisplaybreaks[4]
\usepackage[explicit, pagestyles]{titlesec}
\newcommand{\divisionnameptx}{\relax}%
\newcommand{\authorsptx}{\relax}%
\NewDocumentEnvironment{sectionptx}{mmmmmm}
{%
\renewcommand{\divisionnameptx}{Section}%
\renewcommand{\titleptx}{#1}%
\renewcommand{\subtitleptx}{#2}%
\renewcommand{\shortitleptx}{#3}%
\renewcommand{\authorsptx}{#4}%
\renewcommand{\epigraphptx}{#5}%
\section[#3]{#1}%
\label{#6}%
}{}%
\NewDocumentEnvironment{references-section}{mmmmmm}
{%
\renewcommand{\divisionnameptx}{References}%
\renewcommand{\titleptx}{#1}%
\renewcommand{\subtitleptx}{#2}%
\renewcommand{\shortitleptx}{#3}%
\renewcommand{\authorsptx}{#4}%
\renewcommand{\epigraphptx}{#5}%
\section[#3]{#1}%
\label{#6}%
}{}%
\NewDocumentEnvironment{references-section-numberless}{mmmmmm}
{%
\renewcommand{\divisionnameptx}{References}%
\renewcommand{\titleptx}{#1}%
\renewcommand{\subtitleptx}{#2}%
\renewcommand{\shortitleptx}{#3}%
\renewcommand{\authorsptx}{#4}%
\renewcommand{\epigraphptx}{#5}%
\section*{#1}%
\addcontentsline{toc}{section}{#3}
\label{#6}%
}{}%
\titleformat{\chapter}[display]
{\normalfont\huge\bfseries}{\divisionnameptx\space\thechapter}{20pt}{\Huge#1}
[{\Large\authorsptx}]
\titleformat{name=\chapter,numberless}[display]
{\normalfont\huge\bfseries}{}{0pt}{#1}
[{\Large\authorsptx}]
\titlespacing*{\chapter}{0pt}{50pt}{40pt}
\titleformat{\section}[hang]
{\normalfont\Large\bfseries}{\thesection}{1ex}{#1}
[{\large\authorsptx}]
\titleformat{name=\section,numberless}[block]
{\normalfont\Large\bfseries}{}{0pt}{#1}
[{\large\authorsptx}]
\titlespacing*{\section}{0pt}{3.5ex plus 1ex minus .2ex}{2.3ex plus .2ex}
\titleformat{\subsection}[hang]
{\normalfont\large\bfseries}{\thesubsection}{1ex}{#1}
[{\normalsize\authorsptx}]
\titleformat{name=\subsection,numberless}[block]
{\normalfont\large\bfseries}{}{0pt}{#1}
[{\normalsize\authorsptx}]
\titlespacing*{\subsection}{0pt}{3.25ex plus 1ex minus .2ex}{1.5ex plus .2ex}
\titleformat{\subsubsection}[hang]
{\normalfont\normalsize\bfseries}{\thesubsubsection}{1em}{#1}
[{\small\authorsptx}]
\titleformat{name=\subsubsection,numberless}[block]
{\normalfont\normalsize\bfseries}{}{0pt}{#1}
[{\normalsize\authorsptx}]
\titlespacing*{\subsubsection}{0pt}{3.25ex plus 1ex minus .2ex}{1.5ex plus .2ex}
\usepackage{amsthm}

\tcbset{ theoremstyle/.style={bwminimalstyle, runintitlestyle, blockspacingstyle, after title={\space}, } }
\newtcolorbox[use counter*=cthm]{theorem}[3]{title={{Theorem~\thecthm\notblank{#1#2}{\space}{}\notblank{#1}{\space#1}{}\notblank{#2}{\space(#2)}{}}}, phantomlabel={#3}, breakable, parbox=false, fontupper=\itshape, theoremstyle, }
\tcbset{ corollarystyle/.style={bwminimalstyle, runintitlestyle, blockspacingstyle, after title={\space}, } }
\newtcolorbox[use counter*=cthm]{corollary}[3]{title={{Corollary~\thecthm\notblank{#1#2}{\space}{}\notblank{#1}{\space#1}{}\notblank{#2}{\space(#2)}{}}}, phantomlabel={#3}, breakable, parbox=false, fontupper=\itshape, corollarystyle, }
\NewDocumentEnvironment{introduction}{m}
{\notblank{#1}{\noindent\textbf{#1}\space}{}}{\par\medskip}
\tcbset{ proofstyle/.style={bwminimalstyle, fonttitle=\normalfont\itshape, attach title to upper, after title={\space}, after upper={\space\space\hspace*{\stretch{1}}\(\blacksquare\)}
} }
\newtcolorbox{proofptx}[2]{title={\notblank{#1}{#1}{Proof.}}, phantom={\hypertarget{#2}{}}, breakable, parbox=false, proofstyle }

\usepackage{float}
\usepackage{newfloat}
\usepackage{caption}
\SetupFloatingEnvironment{figure}{fileext=lof,placement={H},within=section,name=Figure}
\makeatletter
\let\c@figure\c@cthm
\makeatother
\usepackage{enumitem}
\newlist{referencelist}{description}{4}
\setlist[referencelist]{leftmargin=!,labelwidth=!,labelsep=0ex,itemsep=1.0ex,topsep=1.0ex,partopsep=0pt,parsep=0pt}
\usepackage[hyperfootnotes=false]{hyperref}
\hypersetup{colorlinks=true,linkcolor=blue,citecolor=blue,filecolor=blue,urlcolor=blue}
\hypersetup{pdftitle={``Pass the Buck'' on a Complete Binary Tree}}

\usepackage{extpfeil}

\title{``Pass the Buck'' on a Complete Binary Tree}
\author{Kenneth Levasseur\\
Department of Mathematical Sciences\\
University of Massachusetts Lowell\\
Lowell, Massachusetts, USA\\
\href{mailto:kenneth_levasseur@uml.edu}{\nolinkurl{kenneth_levasseur@uml.edu}}
}
\date{June 18, 2019}
\begin{document}
\hypertarget{passthebuck}{}
\maketitle
\thispagestyle{empty}
\begin{abstract}
\hypertarget{p-1}{}%
The Stochastic Abacus is employed to compute winning probabilities at each level of the game ``Pass the Buck'' on a complete binary tree with the starting vertex being the root of the tree.  The derivation is also generalized to play on complete \(k\)-ary trees.%
\end{abstract}
\begin{introduction}{Introduction.}%
\hypertarget{p-2}{}%
In the 1970's, Engel \hyperlink{biblio-engel}{[2]} devised the Stochastic Abacus as a way to compute probabilities for certain discrete probability problems with minimal numerical computation.  More recently, Torrence \hyperlink{biblio-torrence}{[6]} used the same technique to determine winning probabilities for players in the game ``Pass the Buck'' for a variety of families of graphs.  The Stochastic Abacus has found more widespread exposure due to a recent article by Propp \hyperlink{biblio-propp}{[4]} in Math Horizons. In this note, the game is analyzed for compete binary trees where the start vertex is the root, and we derive winning probabilities for nodes at different levels. A limiting probability of \(\sqrt{2}-1\) for the root to win as the number of levels goes to infinity is derived. We also observe that the derivation easily generalizes to complete \(k\)-ary trees.%
\end{introduction}%
\typeout{************************************************}
\typeout{Section 1 Pass the Buck}
\typeout{************************************************}
\begin{sectionptx}{Pass the Buck}{}{Pass the Buck}{}{}{section-pass-the-buck}
\hypertarget{p-3}{}%
The game ``Pass the Buck'' is played on a connected undirected graph, with a distinguished ``start vertex.'' The game proceeds in steps starting with the start vertex holding a prize (the ``buck'').  At every stage in the game, the current vertex that holds the buck and its neighboring vertices are selected randomly and uniformly. If the the current vertex is selected, the game ends with that vertex winning. If a neighboring vertex is selected the buck is passed there and process is repeated. More precisely, if the degree of the vertex that holds the buck is \(k\), then the buck moves to any of the neighbors with probability \(\frac{1}{k+1}\) and ends with the player at the current vertex winning with probability \(\frac{1}{k+1}\).%
\end{sectionptx}
\typeout{************************************************}
\typeout{Section 2 The Stochastic Abacus}
\typeout{************************************************}
\begin{sectionptx}{The Stochastic Abacus}{}{The Stochastic Abacus}{}{}{section-stochastic-abacus}
\hypertarget{p-4}{}%
For different graphs, the derivation of the probabilities of any vertex winning can be derived in a variety ways. For example, we can develop a system of equations that determine the probabilities. The game can also be modeled as a Markov Chain and our desired probabilites can be computed using well-know techniques. See Kemeny and Snell \hyperlink{biblio-kemeny}{[3]} for a general introduction to Markov chains and Snell \hyperlink{biblio-snell}{[5]} for a discussion of the connection between Markov chains and the Stochastic Abacus.  Alternatively, the Stochastic Abacus method (also known as Engel's Algorithm) uses only elementary transition rules to compute winning probabilities. We will illustrate all three methods for the case of a complete binary tree up to level 1, \hyperref[fig-binary-1]{Figure~1}, which will serve as a basis for computing the probabilities in larger trees.%
\begin{figure}
\centering
\includegraphics[width=0.6\linewidth]{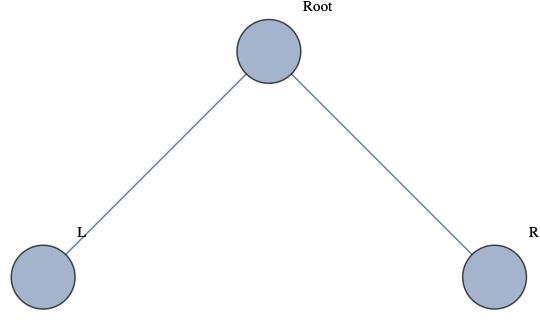}
\caption{A complete binary tree to one level\label{fig-binary-1}}
\end{figure}
\hypertarget{p-5}{}%
In this form of the game, we have three players (root, L, and R) and the buck starts at the root. We can easily derive the probabilities for each vertex by observing that \(p_{\text{root}}= \frac{1}{3}+2\left(\frac{1}{6}p_{\text{root}}\right) \Rightarrow  p_{\text{root}}=\frac{1}{2}\).  By symmetry, \(p_L=p_R=\frac{1}{4}\). Although this is the shortest of our derivations, it doesn't scale so easily.%
\par
\hypertarget{p-6}{}%
Next, we will derive the probabilities using Markov chain theory. The states \(v_x\) corresponding each of the three vertices are states for which the game is in progress, and are non-absorbing states for the process.  We add three absorbing states \(w_x\) to represent winning outcomes for each player. An absorbing state is one for which the process never leaves, representing the end of the game in our case.  The transition matrix for the process with ordering of states \(v_L\), \(v_{\text{root}}\), \(v_R\), \(w_L\), \(w_{\text{root}}\), \(w_R\) is%
\begin{equation*}
T=\left( \begin{array}{cccccc} 0 & \frac{1}{2} & 0 & \frac{1}{2} & 0 & 0 \\ \frac{1}{3} & 0 & \frac{1}{3} & 0 & \frac{1}{3} & 0 \\ 0 & \frac{1}{2} & 0 & 0 & 0 & \frac{1}{2} \\ 0 & 0 & 0 & 1 & 0 & 0 \\ 0 & 0 & 0 & 0 & 1 & 0 \\ 0 & 0 & 0 & 0 & 0 & 1 \\ \end{array} \right)
\end{equation*}
\par
\hypertarget{p-7}{}%
By listing absorbing states last, the form of the transition matrix of an absorbing Markov chain with \(k\) absorbing states is%
\begin{equation*}
\left( \begin{array}{c|c} Q & R \\ \hline \pmb{0} & I_k \\ \end{array} \right)
\end{equation*}
In our  case \(k=3\). If there are a total of \(m\) states, \(Q\) is an \((m-k)\times (m-k)\) matrix of transition probabilities between the non-absorbing states. The transition probabilities from non-absorbing states to absorbing states is contained within \(R\). The matrix \(N R\), where \(N=(I-Q)^{-1}\) is the matrix of probabilities into the different absorbing states. The \(j^{\text{th}}\)row of \(N R\) contains the probabilities of ending in the absorbing states assuming the process starts in state \(j\).%
\par
\hypertarget{p-8}{}%
In our example,%
\begin{equation*}
N R = \left(\left( \begin{array}{ccc} 1 & 0 & 0 \\ 0 & 1 & 0 \\ 0 & 0 & 1 \\ \end{array} \right)-\left( \begin{array}{ccc} 0 & \frac{1}{2} & 0 \\ \frac{1}{3} & 0 & \frac{1}{3} \\ 0 & \frac{1}{2} & 0 \\ \end{array} \right)\right)^{-1}\left( \begin{array}{ccc} \frac{1}{2} & 0 & 0 \\ 0 & \frac{1}{3} & 0 \\ 0 & 0 & \frac{1}{2} \\ \end{array} \right)= \left( \begin{array}{ccc} \frac{5}{8} & \frac{1}{4} & \frac{1}{8} \\ \frac{1}{4} & \frac{1}{2} & \frac{1}{4} \\ \frac{1}{8} & \frac{1}{4} & \frac{5}{8} \\ \end{array} \right)
\end{equation*}
\par
\hypertarget{p-9}{}%
We are mostly concerned with the middle row, which gives probabilities that are consistent with the previous derivation. For this simple case, the graph we have considered is also a path graph and the probabilities were the starting position is either \(L\) or \(R\) is consistent with the general case of ``pass the buck'' on a path graph starting at an end vertex as was considered by Torrence \hyperlink{biblio-torrence}{[6]}.%
\par
\hypertarget{p-10}{}%
The same probabilities we have observed twice will be arrived at using a Stochastic abacus.  The abacus is constructed by first considering the graph to be directed, with each undirected vertex becoming a pair of directed edges. Then we add terminal vertices to the graph, one for each ``internal vertex,'' and an edge leading into each terminal vertex, as in \hyperref[fig-abacus-1]{Figure~2}. The root is designated as the starting position. Initially, we deposit chips into the internal vertices, the number chips being one less than the outdegree of each vertex in the augmented graph, as indicated in each vertex of \hyperref[fig-abacus-1]{Figure~2}.  At this point, the system is ``critically loaded.'' The process then consists of sequentially, adding a chip to the start vertex, \(v_{\text{root}}\), and then repeatedly ``firing'' chips whenever the content of a vertex is greater than or equal to the outdegree of that vertex. This involves distributing a chip along each outgoing edge of the ``loaded vertex'' to neighboring vertices. This process of adding chips and firing as long as possible continues until the chip content of the nonterminal vertices returns to the critical loading state.%
\begin{figure}
\centering
\includegraphics[width=0.6\linewidth]{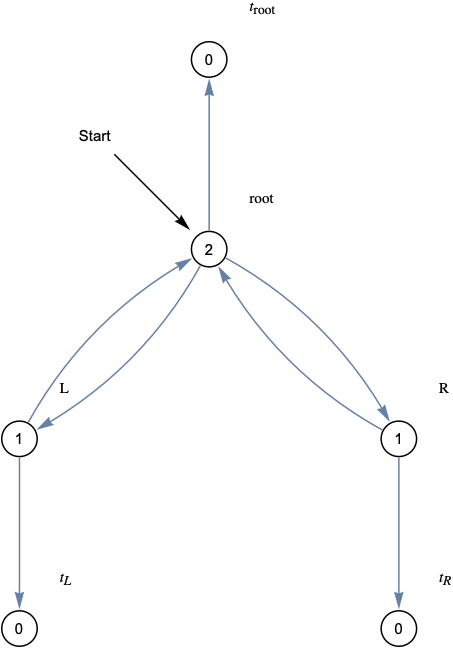}
\caption{Stochastic Abacus for a Complete Level 1 Binary Tree\label{fig-abacus-1}}
\end{figure}
\hypertarget{p-11}{}%
The remarkable fact is that after we have returned to the original critical loading of internal vertices, the probability that any vertex wins the game is equal to the number of chips in its corresponding terminal vertex divided by the total number of chips in all terminal vertices. See Snell \hyperlink{biblio-snell}{[5]} for a proof.  Here is a step by step account of how the process plays out in our example.%
\begin{equation*}
\begin{array}{|c|c|c|c|c|c|c|c|} \hline \text{Step} & \text{Comment} & \text{root} & L & R & t_{\text{root}} & t_L & t_R \\ \hline 1 & \text{Critically } \text{loaded} & 2 & 1 & 1 & 0 & 0 & 0 \\ \hline 2 & \text{Add } 1 \text{ to} \text{root} & 3 & 1 & 1 & 0 & 0 & 0 \\ \hline 3 & \text{root } \text{fires} & 0 & 2 & 2 & 1 & 1 & 1 \\ \hline 4 & L \text{ fires} & 1 & 0 & 2 & 1 & 1 & 0 \\ \hline 5 & R \text{ fires} & 2 & 0 & 0 & 1 & 1 & 1 \\ \hline 6 & \text{Add } 1 \text{ to } \text{root} & 3 & 0 & 0 & 1 & 1 & 1 \\ \hline 7 & \text{root } \text{fires} & 0 & 1 & 1 & 2 & 1 & 1 \\ \hline 8 & \text{add } 1 \text{ to} \text{ root} & 1 & 1 & 1 & 2 & 1 & 1 \\ \hline 9 & \text{add } 1 \text{ to } \text{root} & 2 & 1 & 1 & 2 & 1 & 1 \\ \hline \end{array}
\end{equation*}
It may not be obvious, but if two nodes can fire, as in steps 4 and 5, it doesn't matter in what order they are fired. See  Bjöner \hyperlink{biblio-bjoner}{[1]} for a proof.  After step 9, the three interior vertices are back to being critically loaded, and the process ends.  The total number of chips in the terminal vertices is 4 and the root had 2, so its probability of winning is \(\frac{1}{2}\), consistent with our previous derivation. The other two vertices again have winning probability \(\frac{1}{4}\).%
\end{sectionptx}
\typeout{************************************************}
\typeout{Section 3 Pass the Buck on a Complete Binary Tree}
\typeout{************************************************}
\begin{sectionptx}{Pass the Buck on a Complete Binary Tree}{}{Pass the Buck on a Complete Binary Tree}{}{}{s-ptb-binary-tree}
\hypertarget{p-12}{}%
Consider the game of Pass the Buck on a complete binary tree with \(n\) full levels, \(n\geq 0\), where the buck starts at the root of the tree. We derive formulae for the probabilities that that any node at level \(k\) of the tree, \(0\leq k\leq n\) will win the game. Our derivation is based on an observation that the chips needed at different levels is recursive, with a second order recurrence.%
\begin{theorem}{}{}{theorem-recursion}%
\hypertarget{p-13}{}%
The number of chips in the terminal vertex of the root at the end of the stochastic abacus process, \(a(n)\), follows the recursion \(a(n)=4 a(n-1)-2a(n-2)\), \(n\geq 2\).%
\end{theorem}
\begin{proofptx}{}{proof-1}
\hypertarget{p-14}{}%
If a complete binary tree to \(n\) levels is critically loaded and augmented with terminal vertices, we observe that in order to return the two subtrees starting at level 1 to critical loading status, each vertex at level 1, which are roots of binary trees of level \(n-1\), must fire \(a(n-1)\) times. Every time these vertices fire, they need 4 chips. The initial loading of three chips to each of the vertices at level 1 are used in the first firing, but then to return to critical loading, three other chips are needed. Therefore the root must fire \(4 a(n-1)\) times, almost. There is one other source of chips to each vertex at level 1. That is the vertices at level 2. Each time they fire, they return a chip back up to level 1.  Therefore, the root must fire \(a(n)=4 a(n-1)-2a(n-2)\) times to complete the process. This is then the number of chips that are deposited into the terminal vertex of the root.%
\end{proofptx}
\hypertarget{p-15}{}%
We will apply this recursion together with the two base cases, \(a(0)=1\) and \(a(1)=2\), to compute our probabilities.%
\begin{corollary}{}{}{corollary-1}%
\hypertarget{p-16}{}%
At the end of the stochastic abacus process on a complete binary try of level \(n\), the number of chips deposited into each terminal vertex at level \(k\), \(0\leq k\leq n\) is \(a(n-k)\). The total number of chips that are in all terminal vertices is then \(\sum _{k=0}^n  2^k a(n-k)\).%
\end{corollary}
\begin{corollary}{}{}{corollary-2}%
\hypertarget{p-17}{}%
With the initial conditions \(a(0)=1\) and \(a(1)=2\), \(a(n)=\frac{1}{2} \left(\left(2-\sqrt{2}\right)^n+\left(2+\sqrt{2}\right)^n\right)\), and the total number of chips in all terminal vertices at the end of the stochastic abacus process is \(t(n)=\sum _{k=0}^n  2^k a(n-k)=\frac{\left(2+\sqrt{2}\right)^{n+1}-\left(2-\sqrt{2}\right)^{n+1}}{2 \sqrt{2}}\).%
\end{corollary}
\hypertarget{p-18}{}%
Let \(p(n,k)\) be the probability that one of the \(2^k\) vertices at level \(k\) win the level \(n\) game. The probability that the root is the winner of the level \(n\) game is%
\begin{equation*}
p(n,0)=\frac{a(n)}{t(n)}=\frac{\sqrt{2} \left(\left(2-\sqrt{2}\right)^n+\left(2+\sqrt{2}\right)^n\right)}{\left(2+\sqrt{2}\right)^{n+1}-\left(2-\sqrt{2}\right)^{n+1}}.
\end{equation*}
\par
\hypertarget{p-19}{}%
Interestingly, \(\underset{n\to \infty }{\text{lim}}p(n,n)=\sqrt{2}-1\). More generally, \(p(n,k)=\frac{a(n-k)}{t(n)}\) and \(\underset{n\to \infty }{\text{lim}}p(n,k)=\frac{\sqrt{2}}{\left(2+\sqrt{2}\right)^{k+1}}\).%
\end{sectionptx}
\typeout{************************************************}
\typeout{Section 4 Pass the buck on \(k\)-ary Trees}
\typeout{************************************************}
\begin{sectionptx}{Pass the buck on \(k\)-ary Trees}{}{Pass the buck on \(k\)-ary Trees}{}{}{s-ptb-k-tree}
\hypertarget{p-20}{}%
We can generalize our argument on binary trees to \(k\)-ary trees, \(k\geq 2\). A complete \(k\)-ary tree to level \(n\), \(n>0\), will have a root and \(k\)subtrees, each a \(k\)-ary tree to level \(n-1\). A \(k\)-ary tree up to level 0 is a single vertex, which is the form of a leaf for larger trees. By the same logic as the binary case, if we critically load the augmented directed graph for a complete \(k\)-ary tree up to level \(n\), \(n\geq 2\), then the number of firings of the start vertex that are needed to return to the critically loaded state, \(a(k,n)\), satisfies the recursion \(a(k,n)=(k+2)a(k,n-1)-k a(k,n-2)\). The basis is \(a(k,0)=1\) and \(a(k,1)=2\).%
\end{sectionptx}
\typeout{************************************************}
\typeout{References  References}
\typeout{************************************************}
\begin{references-section-numberless}{References}{}{References}{}{}{references-1}
\begin{referencelist}
\bibitem[1]{biblio-bjoner}\hypertarget{biblio-bjoner}{}Bjöner, A., Lovasz, L., Shor, P. (1991), \textit{Chip-firing games on graphs}, Eur. J. Combin. \textbf{12 (4)}, 283\textendash{}291, doi.org\slash{}10.1016\slash{}s0195-6698(13)80111-4.
\bibitem[2]{biblio-engel}\hypertarget{biblio-engel}{}Arthur Engel (1976), \textit{Why does the probabilistic abacus work?}, Educational Studies in Mathematics \textbf{7}, 59\textendash{}69.
\bibitem[3]{biblio-kemeny}\hypertarget{biblio-kemeny}{}John G. Kemeny and J. Laurie Snell, \textit{Finite Markov Chains}, Undergraduate Texts in Mathematics, Springer- Verlag, New York, 1976.
\bibitem[4]{biblio-propp}\hypertarget{biblio-propp}{}Propp, J. (2018), \textit{Prof.  Engel’s  marvelously  improbable  machines}, Math Horizons, 26(2):  5–9. doi.org\slash{}10.1080\slash{}10724117.2018.1518840.
\bibitem[5]{biblio-snell}\hypertarget{biblio-snell}{}J. Laurie Snell, \textit{The Engel algorithm for absorbing Markov chains}, Available at https:\slash{}\slash{}arxiv.org\slash{}abs\slash{}0904.1413v1
\bibitem[6]{biblio-torrence}\hypertarget{biblio-torrence}{}Bruce Torrence, \textit{Passing the Buck and Firing Fibonacci: Adventures with the Stochastic Abacus}, The American Mathematical Monthly, May 2019, \textbf{126} no.\@\,5, 387\textendash{}399, doi.org\slash{}10.1080\slash{}00029890.2019.1577089.
\end{referencelist}
\end{references-section-numberless}
\end{document}